\documentclass[12pt,a4paper]{article}

\usepackage{amsmath}
\usepackage{amssymb}
\usepackage{t1enc}
\usepackage[latin1]{inputenc}
\usepackage[english]{babel}
\usepackage{amsfonts}
\usepackage{latexsym}
\usepackage{bm}
\newtheorem{theorem}{Theorem}[section]
\newtheorem{prop}[theorem]{Proposition}
\newtheorem{lemma}[theorem]{Lemma}
\newtheorem{cor}[theorem]{Corollary}

\newtheorem{construction}[theorem]{Construction}

\begin{document}
\title{A sharp upper bound for the independence number}

\author{Peter Borg\\[5mm]
Department of Mathematics \\
University of Malta\\
\texttt{p.borg.02@cantab.net}}
\date{} \maketitle

\begin{abstract}
An $r$-graph $G$ is a pair $(V,E)$ such that $V$ is a set and $E$ is a family of $r$-element subsets of $V$. The \emph{independence number} $\alpha(G)$ of $G$ is the size of a largest subset $I$ of $V$ such that no member of $E$ is a subset of $I$. The \emph{transversal number} $\tau(G)$ of $G$ is the size of a smallest subset $T$ of $V$ that intersects each member of $E$. $G$ is said to be \emph{connected} if for every distinct $v$ and $w$ in $V$ there exists a \emph{path} from $v$ to $w$ (that is, a sequence $e_1, \dots, e_p$ of members of $E$ such that $v \in e_1$, $w \in e_p$, and if $p \geq 2$, then for each $i \in \{1, \dots, p-1\}$, $e_i$ intersects $e_{i+1}$). The \emph{degree} of a member $v$ of $V$ is the number of members of $E$ that contain $v$. The maximum of the degrees of the members of $V$ is denoted by $\Delta(G)$. We show that for any $1 \leq k < n$, if $G = (V,E)$ is a connected $r$-graph, $|V| = n$, and $\Delta(G) = k$, then
\[\alpha(G) \leq n - \left \lceil \frac{n-1}{k(r-1)} \right \rceil, \quad \tau(G) \geq \left \lceil \frac{n-1}{k(r-1)} \right \rceil,\]
%
and these bounds are sharp. The two bounds are equivalent.
\end{abstract}

\section{Introduction} \label{Intro}

Unless stated otherwise, we shall use 
small letters such as $x$ to denote non-negative integers or elements (members) of a set. The set $\{1, 2, \dots\}$ of all positive integers is denoted by $\mathbb{N}$. For $1 \leq m \leq n$, we denote $\{i \in \mathbb{N} \colon m \leq i \leq n\}$ by $[m,n]$, and if $m = 1$, then we also write $[n]$. We take $[0]$ to be the \emph{empty set} $\emptyset$.  For a set $X$ and an integer $r \geq 0$, the set $\{Y \subset X \colon |Y| = r\}$ of all $r$-element subsets of $X$ is denoted by ${X \choose r}$. 
Unless stated otherwise, arbitrary sets are assumed to be finite.

A pair $(X,Y)$ is said to be an $r$-graph if $X$ is a set and $Y$ is a subset of $X \choose r$. If $G$ is an $r$-graph $(X,Y)$, then $X$ is represented by $V(G)$ and its members are called \emph{vertices of $G$}, and $Y$ is represented by $E(G)$ and its members are called \emph{edges of $G$}. A $2$-graph is also simply called a \emph{graph}. 

A subset $I$ of $V(G)$ is said to be an \emph{independent set of
$G$} if no edge of $G$ is a subset of $I$. The \emph{independence number of $G$}, denoted by $\alpha(G)$, is the size of a largest independent set of $G$.

A subset $T$ of $V(G)$ is said to be a \emph{transversal of
$G$} if $T$ intersects each edge of $G$ (i.e.~$e \cap T \neq \emptyset$ for each $e \in E(G)$). The \emph{transversal number of $G$}, denoted by $\tau(G)$, is the size of a smallest transversal of $G$. 

Clearly, the complement $V(G)\backslash I$ of an independent set $I$ of $G$ is a transversal of $G$, and the complement $V(G)\backslash T$ of a transversal $T$ of $G$ is an independent set of $G$. By considering the complement of a largest independent set, we obtain $\tau(G) \leq |V(G)| - \alpha(G)$. By considering the complement of a smallest transversal set, we obtain $\alpha(G) \geq |V(G)| - \tau(G)$. Thus, we have 
\begin{equation}\alpha(G) + \tau(G) = |V(G)|. \label{alphatau}
\end{equation}
%


An $r$-graph $H$ is said to be a \emph{subgraph of $G$} if
$V(H)$ is a subset of $V(G)$ and $E(H)$ is a subset of $E(G)$.

A \emph{$vw$-path of $G$} is a sequence $e_1, \dots, e_p$ of edges of $G$ such that $v \in e_1$, $w \in e_p$, and if $p \geq 2$, then for each $i \in [p-1]$, $e_i$ intersects $e_{i+1}$. 
$G$ is said to be \emph{connected} if $G$ has a $vw$-path for every two distinct vertices $v$ and $w$ of $G$. A \emph{component of $G$} is a \emph{maximal} connected subgraph of $G$ (i.e.~a connected subgraph of $G$ that is not a subgraph of another connected subgraph). If $X_1, \dots, X_c$ are pairwise disjoint sets whose union is $X$, then we say that \emph{$X_1, \dots, X_c$ partition $X$}.  It is easy to see that the following holds.

\begin{prop} \label{components} If $G_1, \dots, G_c$ are the distinct components of an $r$-graph $G$, then $V(G_1), \dots, V(G_c)$ partition $V(G)$, and $E(G_1), \dots, E(G_c)$ partition $E(G)$.
\end{prop}


If $v$ and $w$ are distinct vertices in an edge $e$ of $G$, then $v$ and $w$ are said to be \emph{adjacent in $G$}, and we say that $w$ is a \emph{neighbour of $v$ in $G$}, and vice-versa. 
An edge $e$ is said to be \emph{incident to $x$} if $x$ is a member of $e$. For $v \in V(G)$, $N_G(v)$ denotes the set of neighbours of
$v$ in $G$, and the \emph{degree of $v$ in $G$}, denoted by $d_G(v)$, is the number of edges of $G$ incident to $v$. The maximum of the degrees of the vertices of $G$ (i.e.~$\max\{d_G(v) \colon v \in V(G)\}$) is denoted by $\Delta(G)$.

In this paper we provide a sharp (i.e.~attainable) upper bound for the independence number of every connected $r$-graph in terms of the maximum vertex degree. A sharp upper bound for the independence number of every $r$-graph follows immediately. By (\ref{alphatau}), this automatically provides a sharp lower bound for the transversal number.

\section{A sharp upper bound for $\alpha(G)$} \label{upper}

By (\ref{alphatau}), an upper bound for the independence number automatically yields a lower bound for the transversal number, and vice-versa. More precisely, $\alpha(G) \leq |V(G)|-a$ if and only if $\tau(G) \geq a$. Also, $\alpha(G) \geq |V(G)| - a$ if and only if $\tau(G) \leq a$. Various bounds are known for these natural and important parameters. 

A classical theorem of Tur\'{a}n \cite{Turan} says that if $G$ is a $2$-graph and $d$ is the average degree $\frac{1}{|V(G)|}\sum_{v \in V(G)}d_G(v)$, then $\alpha(G) \geq \frac{|V(G)|}{d+1}$. Caro \cite{Caro} and Wei \cite{Wei} independently improved this to $\alpha(G) \geq \sum_{v \in V(G)}\frac{1}{1+d_G(v)}$. Caro and Tuza \cite{CT} generalised the Caro-Wei bound to one for every $r$-graph $G$, given by $\alpha(G) \geq \sum_{v \in V(G)}\prod_{i=1}^{d_G(v)}\left(1 - \frac{1}{(r-1)i+1}\right)$.
Alon \cite{Alon} proved the upper bound $\tau(G) \leq \frac{\ln r}{r}|V(G)| + \frac{|E(G)|}{r}$ for every $r$-graph $G$, and he also showed that the bound is asymptotically sharp; as explained above, this is equivalent to $\alpha(G) \geq \left(1 - \frac{\ln r}{r} \right)|V(G)| - \frac{|E(G)|}{r}$. 

We shall instead prove a sharp upper bound for $\alpha(G)$ in terms of $|V(G)|$ and $\Delta(G)$.

%
%
%
%
For $r \geq 2$ and $k \geq 1$, let $f(r,k)$ be the smallest integer $n$ such that there exists an $r$-graph $G$ with $\Delta(G) = k$ and $|V(G)| = n$. 
\begin{prop}\label{frk} For every $r \geq 2$ and $k \geq 1$, 
\[f(r,k) = \min\left\{m \in \mathbb{N} \colon {m \choose r-1} \geq k\right\} + 1. \]
\end{prop}
\textbf{Proof.} Let $s = \min\left\{m \in \mathbb{N} \colon {m \choose r-1} \geq k\right\}$. Let $G$ be an $r$-graph such that $\Delta(G) = k$ and $|V(G)| = f(r,k)$. Then $d_G(v) = k$ for some $v \in G$. Let $e_1, \dots, e_k$ be the edges of $G$ that are incident to $v$. For each $i \in [k]$, let $e_i' = e_i \backslash \{v\}$. Let $X = \bigcup_{i=1}^k e_i'$. So $e_1, \dots, e_k \in {X \choose r-1}$ and hence $k \leq {|X| \choose r-1}$. So $|X| \geq s$. Since $\{v\} \cup X \subseteq V(G)$ and $v \notin X$, $|V(G)| \geq |X| + 1 \geq s+1$. So $f(r,k) \geq s+1$. Now, since $k \leq {s \choose r-1}$, we can choose $k$ distinct members $a_1, \dots, a_k$ of ${[s] \choose r-1}$. For each $i \in [k]$, let $a_i' = a_i \cup \{s+1\}$. Let $H = ([s+1],\{a_i' \colon i \in [k]\})$. Then $H$ is an $r$-graph with $\Delta(H) = d_G(s+1) = k$ and $|V(G)| = s+1$. So $f(r,k) \leq s+1$. Since $f(r,k) \geq s+1$, we actually have $f(r,k) = s+1$.~\hfill{$\Box$}
\\

In the next section we construct a connected $r$-graph $U_{n,r,k}$ with $\Delta(U_{n,r,k}) = k$, $|V(U_{n,r,k})| = n$ and $\alpha(U_{n,r,k}) = n - \left \lceil \frac{n-1}{k(r-1)} \right
\rceil$ for every $r \geq 2$, $k \geq 2$ and $n \geq f(r,k)$ (see Construction~\ref{const}); we take $U_{r,r,1} = ([r],\{[r]\})$. The following is our main result, which is also proved in the next section.

\begin{theorem}\label{result} Let $r \geq 2$, $k \geq 2$ and $n \geq f(r,k)$. 
If $G$ is a connected $r$-graph such that $|V(G)| = n$ and $\Delta(G) = k$, then
\[\alpha(G) \leq n - \left \lceil \frac{n-1}{k(r-1)} \right
\rceil = \left \lfloor \frac{(k-1)n + 1}{k(r-1)} \right \rfloor,\]
and equality holds if $G = U_{n,r,k}$.
\end{theorem}
%

A graph that consists of only one vertex is called a
\emph{singleton}. For a graph $G$, we denote the set of
non-singleton components of $G$ by $\mathcal{C}(G)$.

\begin{cor}\label{cor1} For every $r$-graph $G$,
\[\alpha(G) \leq |V(G)| - \sum_{H \in \mathcal{C}(G)}\left \lceil
\frac{|V(H)|-1}{\Delta(H)(r-1)} \right \rceil, \]
and equality holds if for each $H \in \mathcal{C}(G)$, $H$ is a
copy of $U_{|V(H)|, r, \Delta(H)}$.
\end{cor}
\textbf{Proof.} Let $s$ be the number of singleton components of
a graph $G$. If $\mathcal{C}(G) = \emptyset$ then $\alpha(G) = s = n$. Suppose $\mathcal{C}(G) \neq \emptyset$. Clearly, $\Delta(H) \geq 1$ for each $H \in \mathcal{C}(G)$. By Theorem~\ref{result}, for any $H \in \mathcal{C}(G)$ with $\Delta(H) \geq 2$ we have $\alpha(H) \leq |V(H)| - \left \lceil \frac{|V(H)|-1}{\Delta(H)(r-1)} \right \rceil$. A connected $r$-graph $K$ with $\Delta(K) = 1$ can only consist of $r$ vertices and an edge containing them (i.e.~$K$ is a copy of $U_{r,r,1}$); thus, for any $H \in \mathcal{C}(G)$ with $\Delta(H) = 1$ we have $\alpha(H) = r-1 = |V(H)| - \left \lceil \frac{|V(H)|-1}{\Delta(H)(r-1)} \right \rceil$. Now, by Proposition~\ref{components}, we clearly have
\begin{align} \alpha(G) &= s + \sum_{H \in \mathcal{C}(G)}
\alpha(H) \leq s + \sum_{H \in \mathcal{C}(G)} \left( |V(H)| -
\left \lceil \frac{|V(H)|-1}{\Delta(H)(r-1)} \right \rceil \right)
\nonumber \\
&= s + \sum_{H \in \mathcal{C}(G)} |V(H)| - \sum_{H \in
\mathcal{C}(G)} \left \lceil \frac{|V(H)|-1}{\Delta(H)(r-1)} \right
\rceil \nonumber \\
&= n - \sum_{H \in \mathcal{C}(G)} \left \lceil
\frac{|V(H)|-1}{\Delta(H)(r-1)} \right \rceil, \nonumber
\end{align}
and by Theorem~\ref{result}, equality holds throughout if each $H
\in \mathcal{C}(G)$ is a copy of $U_{|V(H)|,
\Delta(H)}$.~\hfill{$\Box$}\\

By (\ref{alphatau}), we have the following immediate consequence.

\begin{cor}\label{cor1} For every $r$-graph $G$,
\[\tau(G) \geq \sum_{H \in \mathcal{C}(G)}\left \lceil
\frac{|V(H)|-1}{\Delta(H)(r-1)} \right \rceil, \]
and equality holds if for each $H \in \mathcal{C}(G)$, $H$ is a
copy of $U_{|V(H)|, r, \Delta(H)}$.
\end{cor}
%


\section{Proof of Theorem~\ref{I-bound}}
We start the proof of Theorem~\ref{result} by making the following
observation.
\begin{lemma}\label{I-bound} If $I$ is an independent set of an $r$-graph $G$, then
\[\sum_{v \in V(G) \backslash I} d_G(v) \geq |E(G)|.\]
\end{lemma}
\textbf{Proof.} For each $v \in V(G)$, let $A_v$ be the set of those edges of $G$ that are incident to $v$; so $|A_v| = d_G(v)$. Since $I$ is independent, no edge of $G$ has all its vertices in $I$; in other words, each edge of $G$ has at least one vertex in
$V(G) \backslash I$. So $E(G) = \bigcup_{v \in V(G) \backslash I} A_v$. We therefore have
\[|E(G)| = \left|\bigcup_{v \in V(G) \backslash I} A_v \right| \leq \sum_{v \in V(G) \backslash I} |A_v| = \sum_{v \in
V(G) \backslash I} d_G(v)\]
as required.~\hfill{$\Box$}

\begin{lemma}\label{eremove} If $G$ is a connected $r$-graph, $e \in E(G)$, and $G' = (V(G),E(G) \backslash \{e\})$, then the number of components of $G'$ is at most $r$.
\end{lemma}
\textbf{Proof.} Let $c$ be the number of components of $G'$. Let $G_1, \dots, G_c$ be the components of $G'$. Suppose $e \cap V(G_j) = \emptyset$ for some $j \in [c]$. Then $G_j$ is a component of $G$. Since $G$ is connected, $G = G_j$. So $e \in E(G_j)$ and hence $e \cap V(G_j) = e \neq \emptyset$, a contradiction. So $e \cap V(G_i) \neq \emptyset$ for all $i \in [c]$. Thus, since $V(G_1), \dots, V(G_c)$ partition $V(G)$ (by Proposition~\ref{components}), $|e| = \sum_{i=1}^c |e \cap V(G_i)| \geq \sum_{i=1}^c 1 = c$ and hence $r \geq  c$.~\hfill{$\Box$}

\begin{cor}\label{eremovecor} If $G$ is an $r$-graph, $c$ is the number of components of $G$, $e \in E(G)$, and $G' = (V(G),E(G) \backslash \{e\})$, then the number of components of $G'$ is at most $r + c-1$.
\end{cor}
\textbf{Proof.} Let $G_1, \dots, G_c$ be the components of $G$. By Proposition~\ref{components}, $e \in E(G_i)$ for some $i \in [c]$, and $e \notin E(G_h)$ for each $h \in [c] \backslash \{i\}$. 
Let $G_i' = (V(G_i),E(G_i)\backslash \{e\})$. Then the components of $G'$ are the components of $G_i'$ and the graphs in the set $\{G_h \colon h \in [c] \backslash \{i\}\}$. By Lemma~\ref{eremove}, $G_i'$ has at most $r$ components. So $G'$ has at most $r + c-1$ components.~\hfill{$\Box$}.

\begin{cor}\label{I-bound-cor} If $I$ is an independent set of a
connected $r$-graph $G$, then
\[\sum_{v \in V(G) \backslash I} d_G(v) \geq \frac{|V(G)|-1}{r-1}.\]
\end{cor}
\textbf{Proof.} Let $G$ be a connected $r$-graph, and let $n = |V(G)|$ and $m = |E(G)|$. By Corollary~\ref{eremovecor}, each time an edge is removed from an $r$-graph, the number of components increases by at most $r-1$. Thus, by removing the $m$ edges of $G$ from $G$, the number of components obtained is at most $1 + m(r-1)$; however, the resultant graph is the empty graph $(V(G),\emptyset)$, which has $n$ components (each being a singleton). So $n \leq 1 + m(r-1)$ and hence $m \geq \frac{n-1}{r-1}$. The result
now follows by Lemma~\ref{I-bound}.\\
\\
\textbf{Proof of Theorem~\ref{result}.} Let $G$ be a connected $r$-graph with $|V(G)| = n$ and $\Delta(G) = k$. Let $I$ be a largest
independent set of $G$. By Lemma~\ref{I-bound-cor},
\[\frac{n-1}{r-1} \leq \sum_{v \in V(G) \backslash I} d_G(v) \leq \sum_{v \in
V(G) \backslash I} k = |V(G) \backslash I| k = (n - |I|)k.\]
So $|I| \leq n - \left(\frac{n-1}{k(r-1)} \right)$. Since $|I|$ is an integer and $|I| = \alpha(G)$, we get $\alpha(G) \leq n -
\left\lceil \frac{n-1}{k(r-1)} \right\rceil$ as required.\medskip

We now prove that the upper bound is sharp. Consider the
following construction.

\begin{construction} \label{const} \emph{We define a graph $U_{n,r,k}$ as follows. Let $p = \left\lceil \frac{n-1}{k(r-1)} \right \rceil$. So $n-1 = (p-1)k(r-1)+q$ for some integer $q$ such that $1 \leq q \leq k(r-1)$. Also, $p \geq 1$ since $n > 1$. Let $n' = n-p$. So $n' = (p-1)(k(r-1)-1)+q$. Let $C_{\infty,q,r-1}$ be the Cartesian product $\mathbb{N} \times[k] \times [r-1] = \{(i,j,h) \colon i \in \mathbb{N}, j \in [k], h \in [r-1]$. For each $(i,j,k) \in C_{\infty,q,r-1}$, let
\[x_{i,j,h} = (i-1)(k(r-1) - 1) + (j-1)(r-1) + h.\]
For each $i \in [p]$, let $y_i = n'+i$. For each $(i,j) \in [p] \times [k]$, let $X_{i,j} = \{x_{i,j,h} \colon h \in [r-1]\}$ and $Y_{i,j} = X_{i,j} \cup \{y_i\}$. We have $q = s(r-1)+t$ for some $s \in \{0\} \cup [k]$ and $t \in \{0\} \cup [r-2]$, where $s < k$ if $t > 0$, and $s>0$ if $t=0$. Let $X_{p,s+1}' = [n'-r+2,n']$ and $Y_{p,s+1}' = X_{p,s+1}' \cup \{y_p\}$. Let $M_{n,r,k} = \{Y_{i,j} \colon (i,j) \in ([p-1] \times [k]) \cup (\{p\} \times [s])\} \cup \{Y_{p,s+1}'\}$. If $n > k(r-1) + 1$, then $p \geq 2$ and we take $U_{n,r,k}$ to be $([n],M_{n,r,k})$. Suppose $n \leq k(r-1) + 1$. Then $p=1$, $y_1 = n$, $M_{n,r,k} \subset \{A \in {[n] \choose r} \colon n \in A\}$, $|M_{n,r,k}| \leq s+1 \leq k$, and, since $n \geq f(r,k)$, $k \leq {n-1 \choose r-1}$ by Proposition~\ref{frk}. Thus, we can choose a subset $S_{n,r,k}$ of $\{A \in {[n] \choose r} \colon n \in A\}$ (note that this is a set of size ${n-1 \choose r-1}$) such that $M_{n,r,k} \subseteq S_{n,r,k}$ and $|S_{n,r,k}| = k$, and we take $U_{n,r,k} = ([n],S_{n,r,k})$.}
\end{construction}

We now conclude the proof. Let $U = U_{n,r,k}$ and $M = M_{n,r,k}$. 
We have
\begin{align} x_{1,1,1} &= 1 < \dots < x_{1,1,r-1} = r-1 < x_{1,2,1} = r < \dots < x_{1,k,r-1} = k(r-1) \nonumber \\
&= x_{2,1,1} < \dots < x_{2,k,r-1} = 2k(r-1)-1 \nonumber \\
&= x_{3,1,1} < \dots < x_{3,k,r-1} = 3k(r-1)-2 \nonumber \\
&= x_{4,1,1} < \dots < x_{4,k,r-1} = 4k(r-1)-3 \nonumber \\
& \; \; \vdots \nonumber 
\end{align} 
So $X_{1,1} \cup \dots \cup X_{1,k} = [k(r-1)]$, $X_{2,1} \cup \dots \cup X_{2,k} = [k(r-1),2k(r-1)-1]$, $X_{3,1} \cup \dots \cup X_{3,k} = [2k(r-1)-1,3k(r-1)-2]$, and so on. Let $T = \{y_i \colon i \in [p]\}$ and $J = (\bigcup_{Y \in M} Y) \backslash T$. Then $J = \left(\bigcup_{(i,j) \in [p-1] \times [k]} X_{i,j} \right) \cup \left( \bigcup_{j \in [s]}X_{p,j} \right) \cup X_{p,s+1}' = [n']$ and $T = [n'+1,n]$. Thus, each member of $[n]$ is in some edge of $U$. Now $T$ is a transversal of $U$. If $k \leq k(r-1) + 1$, then $T = \{y_1\}$. If $k > k(r-1) + 1$, then $Y_{1,k}, Y_{2,1}, \dots, Y_{p-1,1}, Y_{p-1,k}, Y_{p,1}$ is a $y_1y_p$ path of $U$ such that each member of $T$ is a member of some edge in the path (note that for each $i \in [p-1]$, $Y_{i,k} \cap Y_{i+1,1} = \{x_{i,k,r-1}\}$ and $Y_{i,1} \cap Y_{i,k} = \{y_i\}$). Thus, since each edge of $U$ is incident to some member of $T$, $U$ is connected. By construction, $U$ is an $r$-graph, $|V(U)| = n$, and $\Delta(U) = d_G(y_1) = k$ (note that $Y_{p,s+1}' = Y_{p,s}$ if $t=0$, and recall that $y_1 = n$ if $n \leq k(r-1)+1$). So $\alpha(U) \leq n - \left\lceil \frac{n-1}{k(r-1)} \right \rceil = n' = |J|$. By (\ref{alphatau}), $J$ is an independent set of $U$ since $J = [n] \backslash T$. So $\alpha(U) \geq |J|$. Since $\alpha(U) \leq |J|$, we actually have $\alpha(U) = |J| = n - \left\lceil \frac{n-1}{k(r-1)} \right \rceil$.~\hfill{$\Box$}

\end{document}